\documentclass[12pt]{article}
\usepackage[english,activeacute]{babel}
\usepackage{amsmath,csquotes,amsfonts,amssymb,amstext,amsthm,amscd,mathrsfs,amsbsy,graphics,graphicx}
\usepackage{xypic}
\usepackage[all]{xy}
\newtheorem{theo}{Theorem}[section]

\newtheorem{con}[theo]{Conjecture}
\theoremstyle{definition}
\newtheorem{defi}[theo]{Definition}

\newtheorem{rem}[theo]{Remark}
\begin{document}
\title{A minimax principle to the injectivity of the Jacobian conjecture}
\author{\\Wei Liu, Quan Xu \footnote {The first author is supported by Ph.D. graduate scholarship in Mathematical Department,Tsinghua University.\
		The second author is supported by post-doc fellowship in Yau Mathematical Sciences center, Tsinghua University.} 
}
\maketitle

\begin{abstract}
The main result of this paper is to prove some type of Real Jacobian Conjecture.  It is proved by the Minimax Principle and asserts if the eigenvalues of $F'(x)$  are bounded from zero and all the eigenvalues of $F'(x)+F'(x)^T$ are strictly same sign, where $ F $ is $ C^1 $ mapping from $ \mathbb{R}^n $ to $ \mathbb{R}^n $, then $ F $ is injective. Moreover $F$ has a $ C^1 $ mapping inverse. 
\end{abstract}

\section{Introduction}
~~~~~ It's well-known that Jacobian Conjecture is first proposed by Keller \cite{kel}.
 \begin{con}{\bf (Jacobian Conjecture)}
	Let $  F(x) $ be  $ k^n \to k^n  $ a polynomial map, where $ k $  is a field of characteristic 0. If the determinant for its jacobian of the polynomial map is a non-zero constant, i.e., $ \det JF(x) \equiv C \in k^{*},~   \forall x\in k^n$.  Then $ F(x) $ has a polynomial inverse map.
\end{con}
  
  For a long study history, it is still open, even for $n=2$. Many results on it, see \cite{BCW}.

  A very important step, for example if $ k=\mathbb{C}^n $, is the following result \cite{CR}.
  \begin{theo}\label{Del,func}
  	Let $F: \mathbb{C}^n\rightarrow \mathbb{C}^n$ is  a polynomial map. 
  	If $ F$ is injective, then $F$ is bijective. Furthermore the inverse is also  a polynomial map.
  \end{theo}  
 The theorem above applied to real case, however, is still open. Generally, the question is transformed to the injectivity of the map in real or complex case. 
Compared with the complex case, we have an analogous conjecture in real case. 
For $k^{n}=\mathbb{R}^n$, the conjecture is
 \begin{con}{\bf(  Real  Jacobian Conjecture)}(RJC)
 	If $F: \mathbb{R}^n\rightarrow \mathbb{R}^n$ is a polynomial map, $\det F'(x)$ is not zero in $\mathbb{R}^n$, then $F$ is a injective map.
 \end{con}
 However, it is a pity that is false and Pinchuk\cite{Pin} constructed a counter-example to (RJC) for $n=2$.

The Pinchuk's counter-example states sufficiently  that
the condition $\det F'(x)$ is not zero in $\mathbb{R}^n$
in (RJC) is much weaker to prove Conjecture 1.1.

In order to enhance 
a sufficient condition for injective, M. Chamberland and 
G. Meisters raise the following conjecture
see (\cite{CM}, Conjecture 2.1):
\begin{con}\label{C1.4}
	Let $F: \mathbb{R}^n\rightarrow \mathbb{R}^n$ be  a $C^1$ map. 
	Suppose there exists an $\epsilon>0$ such that   $ \left| \lambda \right| \ge \varepsilon$ for all the eigenvalues $\lambda$ of $F'(x)$ for all $x\in \mathbb{R}^n$.
	Then F is injective.
\end{con}

\begin{rem}
Pinchuk's counter-example polynomial does not satisfy  the hypotheses of the Conjecture \ref{C1.4}.
\end{rem}
They\cite{CM} also obtain the weak result for the  Conjecture \ref{C1.4}.

\begin{theo}\label{main}
	Let $F: \mathbb{R}^n\rightarrow \mathbb{R}^n$ be  a $C^1$ map. 
	Suppose there exists an $\epsilon>0$ such that   $ \left| \mu  \right| \ge \varepsilon$ for all the eigenvalues $\mu$ of $F'(x)F'(x)^T$ for all $x\in \mathbb{R}^n$.
	Then F is injective.
\end{theo}
In this paper, we prove the conjecture \ref{C1.4} under an additional assumption and the main result is the following theorem.
\begin{theo}\label{main}
	Suppose that $F: \mathbb{R}^n\rightarrow \mathbb{R}^n$ is  a $C^1$ map. 
	If there exists $\epsilon>0$, such that   $ \left| \lambda  \right| \ge \varepsilon ,~ \mu \ge \varepsilon $ or  $\mu \le -\varepsilon$  for all eigenvalues $\lambda,~ \mu$ of $F'(x), F'(x)+F'(x)^T$ respectively, for all $x\in \mathbb{R}^n$.
	 Then F is injective. Moreover F has a $C^1$ inverse map. 
\end{theo}

 In order to prove the theorem \ref{main}, we need to give some definition and notation.

Notation:\\
$\mathbb{C}$: Complex field;\\
$\mathbb{R}$: Real field;\\
$A^T$: the transposition of matrix $A$;\\
$\det A$: the determinant of  matrix $A$;\\
$F'(x)$: The Jacobian matrix of $F(x)$;\\
$X'$: the dual space of $X$;\\ 
$tr A$: the trace of matrix $A$;\\
$\left\| {.} \right\|$: the norm of $\mathbb{R}^n$.

\section{Minimax Principle}

In this section, we will introduce some preparation for the Minimax Principle. 
 
 \begin{defi}($(PS)_c$ sequence)\\
Let $X$ be a Banach space, $J\in C^1(X, \mathbb{R})$.
If ~$\forall \{u_k\}\subset X$, $\exists ~c \in  \mathbb{R} $, such that   \[I({u_k}) \to c,~~I'({u_k}) \to 0, ~\mbox{in}~ X',\]
as $k \to \infty $. Thus $\{u_k\}$ is a $(PS)_c$ sequence of $J$.
\end{defi}

\begin{defi}($(PS)$ condition)\\
	If  ~$\forall (PS)_c$ sequence of $J$ has a convergent subsequence, then $J$ satisfies the $(PS)_c$ condition.
	If $\forall c$, it is said to  satisfy the $(PS)$ condition.
 \end{defi}

\begin{theo}(Minimax Principle)
	Let $X$ be a Banach space, and  $J\in C^1(X, \mathbb{R})$.
	Suppose an open set $ \Omega  \subset X$, $u_0 \in \Omega $,
	and $u_1 \notin  \Omega .$ Set
	\[\Gamma  = \left\{ {\left. {l \in C([0,1],X)} \right|l(i) = {u_i},i = 0,1} \right\}\]
	and
	\begin{equation}
	c = \mathop {\inf }\limits_{l \in \Gamma } \mathop {\sup }\limits_{t \in [0,1]} J(l(t)).
	\end{equation}
	If \\
	(a) $\alpha  = \mathop {\inf }\limits_{\partial \Omega } J(u) > \max \{ J({u_0}),J({u_1})\}$,\\
	(b) $J$ satisfies $(PS)_c$ condition.\\
	Then $c$ is a critical value of $J$.
\end{theo}
In this paper, we use $X=\mathbb{R}^n$ in Theorem 2.3.

\section{The proof of Theorem \ref{main}}
We have already introduced all the necessary ingredients for our theorem.

\begin{proof}
  The proof is by contradiction. Suppose there exist
  $x_1, x_2 \in \mathbb{R}^n, x_1  \ne  x_2  $, such that
  $F(x_1)=F(x_2)$.
  Denote $I\left( x \right) = F\left( {x + {x_1}} \right) - F\left( {{x_2}} \right), \forall x \in \mathbb{R}^n$. So $I(x) \in \mathbb{R}^n$. 
  We define 
  \[J\left( x \right) = I\left( x \right)^T I{\left( x \right)}, ~\forall x \in \mathbb{R}^n.\]
  Thus $J'(x)=2I(x)^TI'(x)$. Let $x_0=x_1-x_2$, then $J(x_0)=J(0)=0$.
 
 Claim： 
 \begin{equation}
 \det I'(x) \ne 0,~\forall x \in \mathbb{R}^n.
 \end{equation}
 
  If there exists $x'\in \mathbb{R}^n$, such that $detI'(x') = 0$.
  So $detF'(x'+x_1)=0.$
  Thus, $\exists~ y \in \mathbb{R}^n, y\ne 0 $, s.t.
  \begin{equation}
   F'(x'+x_1)y=0.
  \end{equation}
  That is $F'(x'+x_1)$ has a zero eigenvalue. It contradicts the eigenvalues of $F'(x)$ are bounded from zero. 
  
  Thus $I'(x)$ is a invertible matrix.
  If $J'(x) = 0,~\forall x \in \mathbb{R}^n $, 
  i.e. $2I(x)^TI'(x)=0,~\forall x \in \mathbb{R}^n $, thus $I'(x)I(x)=0,~\forall x \in \mathbb{R}^n.$
  $I(x)=0$. So $J(x)=0$.

  Next, we prove $J$ satisfies the geometric condition
  (a) in Theorem 2.3.
  
  Since $J(0)=J(x_0)=0$, it is sufficient to prove 
  $\exists ~ r$, such that 
  \begin{equation}
   J(u) > 0,\forall u \in \partial {B_r}(0).
  \end{equation} 
   Claim:
  $x=0$ is an isolated zero point of $J(x)$.
  
  For each component $I_i(x)$ of $I(x)$, so 
  $I_i(x)=I'(y_i)x$, here $y_i$ connects $0$ to $x$, $i=1,2...n$.
  Define a continuous function $\beta (x)$ as
  \[\beta \left( x \right) = \left\{ \begin{array}{l}
  {({{I'}_1}({y_1}),{{I'}_2}({y_2})...{{I'}_n}({y_n}))^T},x \ne 0,\\
  I'(0),~~~~~~~~~~~~~~~~~~~~~~~~~~~~x = 0.
  \end{array} \right.\]
  Thus $I(x)=\beta (x)x$, $\forall x\in \mathbb{R}^n$.
  Define 
  \[\gamma ({x_1},{x_2}...,{x_n}) = {({{I'}_1}({x_1}),{{I'}_2}({x_2})...,{{I'}_n}({x_n}))^T}.\]
  Thus $\gamma(x,x...x)=I'(x)$ and $\gamma(y_1,y_2...,y_n)=\beta(x)$.
  Therefore \[\det \gamma (0,0,...,0) = \det I'(0) \ne 0.\]
  By the continuity of $\gamma$, there exists a positive number $r>0$, such that
  \[\det \gamma ({x_1},{x_2}...,{x_n}) \ne 0,~ \mbox{for}~ (x_1,x_2,...,x_n)\in B_r(0). \]
  Thus $det \beta (x)\ne 0, ~\forall x\in B_{r/{\sqrt n }}(0).$
  Therefore $0$ is a isolated zero point of $I(x)$.
  
  Let $\alpha  = \mathop {\inf }\limits_{\partial B_{r/\sqrt{n}}(0) } J(x)$.
  It is a positive number since $J(x)$ is continuous and is not zero on $\partial B_{r/\sqrt{n}}(0).$
  
  Thus $J(x)$ satisfies the condition (a) in theorem 2.3.
  
  If $c$ is a critical value of $J$, that is $\exists~ x_c\in \mathbb{R}^n $， such that $J'(x_c)=0$. Thus
  \[0 < \alpha  \le c = J({x_c}) = 0.\]
  Obviously, it's impossible. 
  
  By Theorem 2.3, the condition (b) i.e. $(PS)_c$ condition
  does not hold.
  There is  a sequence $\{x_k\} \subset \mathbb{R}^n $, such that
 \[ (i)J({x_k}) \to c;~~
  (ii)J'({x_k}) \to 0;~~
  (iii)\left\| {x_k} \right\| + \infty.\]
  
  If   $\mu \le -\varepsilon$  for all eigenvalues $\mu$ of $F'(x)+F'(x)^T$  for all $x\in \mathbb{R}^n$.
  
 Let  $\mu_0$ denote the maxmum eigenvalue of a Hermitian matrix $A$. Define
 \begin{equation}\label{key1}
{\mu _0} = \mathop {\sup }\limits_{Y \ne 0} \frac{{{Y^T}AY}}{{{Y^T}Y}}.
 \end{equation}
Set $A = I'(x_k) + I'(x_k)^T ~\mbox{and}~ Y=I(x_k)$. By (\ref{key1}), we obtain 
\begin{equation}\label{key2}
\begin{aligned}
{\mu _0}\left( {{x_k}} \right) &\ge \frac{{I{{({x_k})}^T}(I'({x_k}) + I'{{({x_k})}^T})I({x_k})}}{{I{{({x_k})}^T}I({x_k})}}\\
& = \frac{{I{{({x_k})}^T}I'({x_k})I({x_k}){\rm{ + }}I{{({x_k})}^T}I'{{({x_k})}^T}I({x_k})}}{{I{{({x_k})}^T}I({x_k})}}
\\
& = \frac{{2I{{({x_k})}^T}I'({x_k})I({x_k})}}{{I{{({x_k})}^T}I({x_k})}}.
\end{aligned}
\end{equation}
By (i), one gets
\begin{equation}\label{key3}
 I{({x_k})^T}I({x_k}) = J({x_k}) \to c > 0.
\end{equation}
 By (ii) and (\ref{key3}), we obtain
 \begin{equation}\label{key4}
   \begin{aligned}
 2I{({x_k})^T}I'({x_k})I({x_k}) &\le 2\left\| {I{{({x_k})}^T}I'({x_k})} \right\|\left\| {I({x_k})} \right\|\\
 &=\left\| {J'({x_k})} \right\|{(I{({x_k})^T}I({x_k}))^{\frac{1}{2}}}\\
 & = \left\| {J'({x_k})} \right\|\sqrt {J({x_k})}\to 0.
 \end{aligned}
 \end{equation}
Combining (\ref{key2}), (\ref{key3}) with (\ref{key4}), as $k\to +\infty $,
 thus 
 \begin{equation}\label{key5}
  \mu_0(x_k) \ge 0.
 \end{equation}
 By $\mu \le -\varepsilon$ for all the eigenvalue $\mu $ of  $F'(x)+F'(x)^T$, thus 
  \begin{equation}\label{key6}
   {\mu _0}\left( {{x_k}} \right)<0.
  \end{equation}
 In (\ref{key5}) and (\ref{key6}), letting $k\to +\infty $,
  thus,
  \[{\mu _0}\left( {{x_k}} \right) \to 0.\]
  It contradicts.
  
   If   $\mu \ge \varepsilon$  for all eigenvalues $\mu$ of $F'(x)+F'(x)^T$  for all $x\in \mathbb{R}^n$.
 
   Let $ {\mu _1} = \mathop {\inf }\limits_{Y \ne 0} \frac{{{Y^T}AY}}{{{Y^T}Y}}.$
   By the same method, we obtain contradiction.
    
   Therefore $F$ is injective.
  Thus $F$ has a inverse map, denoted by $G$. Since
  $F \cdot G = \text{Id}$ and $F'$ exists, then if $\forall x \in \mathbb{R}^n$, such that $ F \cdot G (x)=x$, then
  \[F'(G(x)) \cdot G'(x) = 1.\]
  Since $F'$ exists and $\det F'\ne 0$, $G'(x)$ exists, and
  \[G'(x) = F'{(G(x))^{ - 1}}.\]
  By $F\in C^1$, we obtain $G\in C^1$.
  
\end{proof}
\begin{rem}
	A further consideration to the question, it is either to give an analogue to minimax principle in complex case by our theorem or to give an real analogue to theorem 1.2, which is great progress to the 
	Jacobian conjecture.
\end{rem}

\section*{Acknowledgements}
The results in this paper were obtained  from some discussion. I want to express my sincere gratitude to Professor M. Chamberland and G. Meisters. Also, we are grateful to  Professors A. van den Essen and D. Wright, who  provided helpful comments. Finally, I would like to thank any readers who can give amount of precious opinions for improvement of this paper.

Address of authors:\\

   Wei Liu 
\\ A322, Department of Mathematical Sciences, Tsinghua university,
 100084, Haidian District, Beijing China.
\\Email: liuw16@mails.tsinghua.edu.cn.\\

  Quan Xu
\\ Jingzhai 204, Yau mathematical Sciences center, Tsinghua University, 
100084, Haidian District, Beijing China.
\\Email: xqmathxy@mail.tsinghua.edu.cn.


\newpage\begin{thebibliography}{XXX}
\bibitem[BCW82]{BCW}  H. Bass, E. Connell and D. Wright. \textit{The Jacobian Conjecture: reduction of degree and formal expansion of the inverse.} Bull. Amer. Math. Soc. 7 (1982), no. 2, 287-330.

\bibitem[BV05]{BV}M. de Bondt and A. van den Essen.\textit{ A Reduction of the Jacobian Conjecture to the Symmetric Case.}  Proc. Amer. Math. Soc. 133 (2005), no. 8, 2201-2205.
	
\bibitem[CM98]{CM} M. Chamberland and G. Meisters. \textit{A mountain pass to the Jacobian Conjecture.} Canad. Math. Bull. Vol.41 (1998), no 4, 442-451.

\bibitem[CR91]{CR} Cynk and K. Rusek. \textit{Injective endomorphisms of algebraic and analytic sets.} Ann. Polonici Math. 56, no 1, (1991). 

\bibitem[Ess97]{Ess} A. van den Essen. \textit{To believe or not to believe: the Jacobian Conjecture.} 
Rend. Sem. Mat. 55 (1997), no 4.

\bibitem[Fri99]{Fri}  S. Friedland . \textit{Monodromy, diﬀerential equations and the Jacobian conjecture.} 
Ann. Polonici Math. 72, no 3, (1999).

\bibitem[Kel39]{kel}O. H. Keller.\textit{ Ganze Gremona-Transformation.}   Monats. Math. Physik 47 (1939), 299-306.
 
\bibitem[Pin94]{Pin} S. Pinchuk. \textit{A counterexamle to the strong real Jacobian conjecture.} 
Math. Zeitschrift. 217 (1994), 1-4.


\bibitem[Zha07]{Zha}  W. Zhao. \textit{Hessian nilpotent polynomials and the Jacobian conjecture.} Tran. Amer. Math. Soc. 359 (2007), no. 1, 249-274.

\end{thebibliography}
\end{document}